\documentclass[12pt,reqno]{amsart}
\usepackage{amsmath,amssymb,amsthm, xcolor}
\usepackage{geometry}
\newtheorem{theorem}{Theorem}[section]
\newtheorem{lemma}[theorem]{Lemma}

\newtheorem*{thmA}{Theorem A}

\title[Indices divisible by at most two primes]{Finite simple groups whose maximal subgroups have indices divisible by at most two distinct primes}

    \author[Beltr\'an and Shao]{Antonio Beltr\'an\\
    Departamento de Matem\'aticas\\
    Universitat Jaume I \\
  12071 Castell\'on\\
   Spain\\ 
    \\Changguo Shao\\
    College of science\\
NANJING UNIVERSITY OF POSTS AND TELECOMMUNICATIONS\\
NANJING 210023 YADONG
        \\
      China}



\date{}

\begin{document}

\begin{abstract}
A theorem by Guralnick classifying finite simple groups with subgroups of prime power index implies that $\mathrm{PSL}(2,7)$ is the unique finite non-abelian simple group whose maximal subgroups all have prime power index. Motivated by this result and appealing to Li and Li's classification of primitive permutation groups of degree $p^a q^b$, where $p$ and $q$ are distinct primes, we classify all finite non-abelian simple groups whose maximal subgroups have indices divisible by at most two distinct primes, proving  that they are precisely $A_5$, $A_6$, $A_7$, $A_8$, ${\rm PSL}(2,7)$, ${\rm PSL}(2,8)$, ${\rm PSL}(2,11)$, ${\rm PSp}(4,3)$, ${\rm PSU}(3,3)$, and ${\rm PSU}(3,4)$.
\end{abstract}

\maketitle

\section{Introduction}

The relationship between the arithmetic properties of the indices of maximal subgroups and the structure of a finite group is a classical theme in finite group theory. One of the earliest results in this area is attributed to P. Hall, who proved that a finite group $G$ is solvable if the index of every maximal subgroup of $G$ is either a prime or the square of a prime (for instance, \cite[9.3.1]{Robinson}).  While  all maximal subgroups in a solvable group  have  prime power indices, the  exponents need not be restricted to $1$ or $2$ as in P. Hall's theorem; indeed, maximal subgroups of index $p^n$ for any prime $p$ and any positive integer $n$ may exist in a solvable group.  This occurs in every solvable Frobenius group having a kernel of order $p^n$ and a complement acting irreducibly on the kernel; for instance, the affine group ${\rm AGL}(1,p^n)$.

Beyond solvable groups, a natural problem is to classify finite simple groups whose maximal subgroups satisfy certain arithmetic restrictions on their indices. In 1983, Guralnick \cite{Gur} classified all  finite non-abelian simple groups possessing a subgroup of prime power index. As a result, he determined all finite non-abelian simple groups in which every maximal subgroup has prime power index: up to isomorphism, the only such group is  $\mathrm{PSL}(2,7) \cong \mathrm{PSL}(3,2)$ (whose maximal subgroups have indices $7$ and $8$).  Furthermore, he characterized the structure of a finite group $G$ in which every maximal subgroup has prime power index, showing that $G$ is either solvable or satisfies $G/S(G) \cong \text{PSL}(2,7)$, where $S(G)$ denotes the solvable radical  of $G$.

Extending these ideas from one prime to two primes, it is natural to ask for a classification of finite non-abelian simple groups $S$ such that the index $|S:M|$ of every maximal subgroup $M$ of $S$ is divisible by at most two distinct primes (depending on $M$). 
This study is closely related to the classification of primitive permutation groups of degree $p^a q^b$, where $p$ and $q$ are primes. This was achieved by C.H. Li and X.H. Li \cite{LL}, who also provided a complete classification of all maximal subgroups of finite simple groups having index $p^a q^b$.   Relying on their work, in this paper we give a complete classification of all finite non-abelian simple groups whose maximal subgroups all have indices divisible by at most two distinct primes.

\begin{thmA}\label{main}
Let $S$ be a finite non-abelian simple group. Every maximal subgroup of $S$ has index divisible by at most two distinct primes if and only if $S$ is isomorphic to one of the following groups:
\begin{enumerate}
    \item[\textnormal{(1)}] $A_5$, $A_6$, $A_7$, $A_8$;
    \item[\textnormal{(2)}] ${\rm PSL}(2,7)$, ${\rm PSL}(2,8)$, ${\rm PSL}(2,11)$;
    \item[\textnormal{(3)}] ${\rm PSp}(4,3)$;
   \item[\textnormal{(4)}] ${\rm PSU}(3,3)$ or ${\rm PSU}(3,4)$.
\end{enumerate}
\end{thmA}

In contrast to the prime-power case, however, Theorem A does not directly yield a description of finite non-solvable groups whose maximal subgroups have indices divisible by at most two primes. The main difficulty is that this arithmetic condition does not force all non-abelian composition factors of the group to belong to the list in Theorem A.
For instance, consider the almost simple group $G = \operatorname{PGL}(2, 17)$. Its maximal subgroups have indices $2$ (corresponding to ${\rm PSL}(2,17))$, $18 = 2 \cdot 3^2$, $136 = 2^3 \cdot 17$, and $153 = 3^2 \cdot 17$, each of which is divisible by at most two distinct primes, whereas its socle, $\operatorname{PSL}(2, 17)$, has a maximal subgroup of index $102 = 2\cdot 3 \cdot 17$.  Thus, $G$ satisfies the condition, but has a non-abelian composition factor that falls outside the list of Theorem A. In particular, this arithmetic condition is not inherited by normal subgroups, so extending Theorem A to non-solvable groups remains open.
\section{Proof}

The proof of Theorem A relies on analyzing the indices of maximal subgroups in finite non-abelian simple groups, following the work of Li and Li \cite{LL} and  the classifications and descriptions given in \cite{BHR,KL}. The analysis is naturally divided into three cases according to the Classification of Finite Simple Groups: 
 alternating groups (Theorem \ref{thm2.1}), sporadic groups (Theorem \ref{thm2.2}), and groups of Lie type (Theorem \ref{thm2.4}). If $n$ is a positive integer, we denote by $\pi(n)$ the set of primes dividing $n$.

\begin{theorem}\label{thm2.1}
Let $G = A_n$ be the alternating group of degree $n \ge 5$. If $|\pi(|G:M|)|\leq 2$ for every maximal subgroup $M$ of $G$, then $G$ is one of the following groups: $A_5, A_6, A_7$, or $A_8$.
\end{theorem}

\begin{proof}

By \cite{Con}, it is straightforward to check that $A_n$ satisfies the hypothesis for $n \in \{5, 6, 7, 8\}$. Suppose by way of contradiction that $n \geq 9$ and that $|\pi(|G:M|)| \leq 2$ for every maximal subgroup $M$ of $G$.

Since $n \geq 9$, the stabilizer of a $3$-element subset, $H = (S_3 \times S_{n-3}) \cap A_n$, is a maximal subgroup of $G$ of index $|G:H|  = \binom{n}{3}$ (see, for instance, \cite{LPS87}). Because this index is larger than $n$, Guralnick's classification \cite[Theorem 1]{Gur} implies that  $\binom{n}{3}$ is not a prime power. Alternatively,  this fact can be checked by a standard  elementary argument that we omit. Thus, our hypothesis forces $|\pi(|G:H|)| = 2$. However, $H$ does not appear in the classification of subgroups of $A_n$ ($n \geq 9$) whose index has two prime divisors given in \cite[Lemma 3.1 and Table 3.1]{LL} (in fact, Case 1 in the proof of \cite[Lemma 3.1]{LL} shows that $\binom{n}{3}\neq p^aq^b$). This contradiction completes the proof.
\end{proof}

\begin{theorem}\label{thm2.2} There are no sporadic simple groups all of whose maximal subgroups have indices divisible by at most two distinct primes.
\end{theorem}

\begin{proof}
Assume by way of contradiction that $G$ is a sporadic simple group satisfying $|\pi(|G:M|)| \leq 2$ for every maximal subgroup $M$ of $G$. According to \cite[Table 3.2]{LL}, $G$ must be isomorphic to one of the Mathieu groups $M_{11}, M_{12}, M_{22}, M_{23}, M_{24}$, or to $J_2, Co_2, HS$, or $McL$. We reach a contradiction in each case by exhibiting a maximal subgroup $M$ from \cite{Con} with $|\pi(|G:M|)| = 3$.

 If $G$ is a Mathieu group, then $M$ can be chosen to be $S_5$ in $M_{11}$ or $M_{10}:2 \cong A_6 \cdot 2^2$ in $M_{12}$, both of index $2 \cdot 3 \cdot 11$; $2^4:S_5$ in $M_{22}$ of index $3 \cdot 7 \cdot 11$;  $A_8$ in $M_{23}$ of index $2 \cdot 11 \cdot 23$; and  $M_{22}:2$ in $M_{24}$ of index $2^2 \cdot 3 \cdot 23$.
    
For the remaining four candidates, we may take $M \cong 3 \cdot \mathrm{PGL}(2,9)$ in $J_2$ (index $2^3 \cdot 5 \cdot 7$), $M \cong \mathrm{PSU}(6,2):2$ in $Co_2$ (index $2^2 \cdot 5^2 \cdot 23$), $M \cong S_8$ in $HS$ (index $2^2 \cdot 5^2 \cdot 11$), and $M \cong \mathrm{PSU}(3,5)$ in $McL$ (index $2^3 \cdot 3^4 \cdot 11$).

Therefore, all sporadic simple groups are ruled out.
\end{proof}

Before dealing with simple groups of Lie type,
recall that if $b$, $n$ are integers greater than $1$, then a prime $p$ is called a primitive prime divisor or a Zsigmondy
prime for $(b,n)$ if $p \nmid b$ and the order of $b$ $({\rm mod}$ $p)$ equals $n$. This means that $p \mid b^n-1$ but $p$ does not divide $b^i-1$ for every $1\leq i \leq n-1$. We will repeatedly use the following well-known properties.

\begin{lemma}{\cite{Feit}}  \label{ppd} Let $b$ and $n$ be positive integers greater than $1$.

\begin{itemize}
\item[(1)] There are primitive prime divisors of $b^n-1$ unless $(b,n) = (2,6)$ or $b=2^r-1$, with $r\geq 2$ and $n = 2$.

\item[(2)]  Each primitive prime divisor $p$ of $b^n-1$ is at least $n + 1$. 
\end{itemize}
\end{lemma}

\begin{theorem}\label{thm2.4}
Let $G$ be a finite non-abelian simple group of Lie type. Every maximal subgroup of $G$ has index divisible by at most two distinct primes if and only if $G$ is isomorphic to one of the following groups:
\begin{enumerate}
   
    \item[\textnormal{(1)}] ${\rm PSL}(2,5)$, ${\rm PSL}(2,7)$, ${\rm PSL}(2,8)$, ${\rm PSL}(2,9)$, ${\rm PSL}(2,11)$,  ${\rm PSL}(4,2)$;
    \item[\textnormal{(2)}] ${\rm PSp}(4,3)$;
   \item[\textnormal{(3)}] ${\rm PSU}(3,3)$ or ${\rm PSU}(3,4)$.
\end{enumerate}
\end{theorem}

\begin{proof}

 If $|\pi(|G:M|)|=1$ for every maximal subgroup $M$ of $G$, then $G\cong {\rm PSL}(2,7)$ by \cite[Corollary 3]{Gur}.  Henceforth, we  assume that $G$ has at least one maximal subgroup with index divisible by two distinct primes.  If $K$ is a maximal subgroup  of $G$ with $|\pi(|G:K|)|=2$, then the pair $(G,K)$ must appear in the classification given in \cite{LL}. More precisely, the candidates for $(G,K)$ are listed in Tables 4.1, 4.2, 4.3, 4.5, and 5.1 of \cite[Theorem 1.1]{LL}, corresponding to linear, symplectic, unitary, orthogonal, and exceptional groups, respectively. We analyze each of these five cases.

\medskip 

(1) $G$ is a linear group, that is,  ${\rm PSL}(n,q)$, where $q=s^f$ and $s$ is a prime. Let $Z=Z({\rm SL}(n,q))$ and $d=|Z|$. We analyze the subcases $n=2$, $n=3$, and $n\geq 4$ separately.

\medskip
(1.1) $n=2$, that is, $G\cong {\rm PSL}(2,q)$.
If $q\leq 11$, then  we have $q=4,  5, 7,  8, 9$ or $11$ and the groups ${\rm PSL}(2,4)\cong{\rm PSL}(2,5)$, ${\rm PSL}(2,7)$, ${\rm PSL}(2,8)$, ${\rm PSL}(2,9)$, and ${\rm PSL}(2,11)$ satisfy the hypothesis of the theorem (by \cite{Con}).  Thus, we only consider $q>11$ in what follows.

\medskip
Suppose first that $q$ is even, that is, $q=2^f$ with $f\geq 4$. We know from \cite[Table 8.1]{BHR} that $G$ has maximal subgroups isomorphic to $D_{2(q-1)}$ and $D_{2(q+1)}$, whose indices are $q(q+1)/2$ and $q(q-1)/2$, respectively. Since $q/2=2^{f-1}\geq 8$ and $(q,q\pm 1)=1$, the hypothesis on the indices of maximal subgroups implies that $2^f+1=p_1^{a_1}$ and $2^f-1=p_2^{a_2}$, where $p_1$ and $p_2$ are primes and $a_1$ and $a_2$ are positive integers. Assume first that $a_2=1$. Then $p_2=2^f-1$ is a Mersenne prime, and this forces $f$ to be an odd prime (as $f\geq 4$, we have $f\geq 5$). Accordingly, $p_1^{a_1}=2^f+1$ is divisible by $3$, so $p_1=3$ and $3^{a_1}-1=2^f$. Since $f\geq 5$, we have $a_1\geq 4$, and Lemma \ref{ppd} ensures that $3^{a_1}-1$ has a primitive prime divisor $r$ with $r\geq a_1+1\geq 5$. The divisibility $r\mid 2^f$ forces $r=2$, a contradiction. Therefore, we may assume $a_2>1$. Then $p_2^{2a_2}-1$ has a primitive prime divisor $r_2$ with $r_2\geq 2a_2+1\geq 5$. Using $p_2^{2a_2}-1=(p_2^{a_2}+1)(p_2^{a_2}-1)$ together with $r_2\nmid p_2^{a_2}-1$, we get $r_2\mid p_2^{a_2}+1=2^f$, so $r_2=2$, contradicting $r_2\geq 5$.

\medskip
We now suppose that $q$ is odd.
Assume first that $f=1$. As above,  $G$ has maximal subgroups $D_{q-1}$ and $D_{q+1}$ with indices $q(q+1)/2$ and $q(q-1)/2$, respectively. Given that $q$ is coprime to $q+1$ and $q-1$, by hypothesis, both $(q-1)/2$ and $(q+1)/2$ are prime powers. Moreover, as $4\mid (q^2-1)$, either $(q-1)/2$ or $(q+1)/2$ is even. Consider first the case where $(q-1)/2$ is even. Then $q-1=s-1=2^{b_1+1}$, where $b_1$ is a positive integer (in fact, $b_1 \geq  3$  as $ q > 11)$, and $(s+1)/2=s_2^{b_2}$, where $s_2$ is a prime and $b_2$ is a positive integer. This yields $s=2^{b_1+1}+1$, which in turn implies that $s$ is a Fermat prime and $b_1+1=2^t$, where $t$ is a positive integer. On the other hand, $s=2\cdot s_2^{b_2}-1$, so we obtain  $2\cdot s_2^{b_2}-1=2^{b_1+1}+1$. Consequently, $s_2^{b_2}=2^{b_1}+1$ (so $s_2$ is odd). If $b_2=1$, then $s_2=2^{b_1}+1$ is a prime, from which we deduce that $b_1$ is a power of $2$, contradicting the fact that $b_1 +1=2^t$. If $b_2=2$, then $s_2^2-1=(s_2+1)(s_2-1)=2^{b_1}$, and this implies that $(s_2-1)\mid (s_2+1)=(s_2-1)+2$. It follows that $s_2=3$, so $q=s=17$. However, by \cite{Con}, ${\rm PSL}(2,17)$ has a maximal subgroup isomorphic to $S_4$ with index $102=2\cdot3\cdot17$,  a contradiction. If $b_2 \geq 3$, then by Lemma \ref{ppd}, $s_2^{b_2} - 1$ admits a primitive prime divisor $r \geq b_2 + 1 \geq 4$ dividing $2^{b_1}$, again a contradiction. This finishes the case where $(q-1)/2$ is even. 
Assume now that $(q+1)/2$ is even and write $(q+1)/2=2^{l_1}$, where $l_1$ is a positive integer. Hence $q =s =2^{l_1+1}-1$ is a prime. In particular, $l_1 +1$ is prime as well. In addition, as $q=s>11$, we have $l_1+1\geq 5$,  so $l_1\geq 4$ is even.  Write $l_1=2m$, where $m\geq2$ is an integer. Since $(q-1)/2$ is also a prime power, we write $(q-1)/2=p_1^{l_2}$, where $p_1$ is a prime and $l_2$ is a positive integer, that is, $q=s=2p_1^{l_2}+1$. Thus, we have $2p_1^{l_2}+1=2^{l_1+1}-1$, that is, $p_1^{l_2}=2^{l_1}-1$.  The fact that $l_1=2m$ yields $p_1^{l_2}=2^{l_1}-1=(2^m+1)(2^m-1)$. Hence we obtain $(2^m-1)\mid (2^m+1)=(2^m-1)+2$. As a consequence, $(2^m-1)\mid 2$, whence $m=1$, contradicting $m\geq 2$.

\medskip
Suppose next that $f>1$.  
Since $q>11$  is odd, $G$ has maximal subgroups $D_{q-1}$ and $D_{q+1}$ by \cite[Table 8.1]{BHR} again. Furthermore, $|G:D_{q-1}|=q(q+1)/2 $ and $|G:D_{q+1}|=q(q-1)/2 $. Since $q$ is coprime to $q+1$ and $q-1$, both $(q+1)/2$ and $(q-1)/2$ are prime powers. Write $(q-1)/2=s_1^{b_1}$ and $(q+1)/2=s_2^{b_2}$, where $s_1$ and $s_2$ are primes, and $b_1$ and $b_2$ are positive integers. Since  $q$ is odd, $4\mid (q^2-1)$, so either $(q+1)/2$ or $(q-1)/2$ is divisible by 2, which yields either $s_1=2$ or $s_2=2$. If $s_1=2$, then $q-1=2s_1^{b_1}$, that is, $s^f-1=2^{b_1+1}$.  If $f=2$, then $(s+1)(s-1)=2^{b_1+1}$, which means that both $s+1$ and $s-1$ are powers of $2$. It follows that $s-1\mid s+1$, which forces $s=3$ and $q=9$, a contradiction. If $f>2$, then, by Lemma \ref{ppd}, we have that $s^f-1$ has a primitive prime divisor $r_0$ with $r_0\geq f+1\geq4$. Consequently, $r_0\mid (s^f-1)=2^{b_1+1}$, a contradiction. 
 Finally, if $s_2=2$, then $s^f+1=2^{b_2+1}$.   Since $s$ is odd and $2f\geq4$, we know from Lemma \ref{ppd}  that $s^{2f}-1$ has a primitive prime divisor $r_1$ with $r_1\geq 2f+1$. From $r_1\mid (s^{2f}-1)=(s^f+1)(s^f-1)$ and $r_1\nmid (s^f-1)$, we conclude that $r_1\mid (s^f+1)=2^{b_2+1}$, also a contradiction.

\medskip
(1.2) $n=3$, that is, $G\cong \mathrm{PSL}(3,q)$. If $q=2$, then $G\cong \mathrm{PSL}(3,2)\cong \mathrm{PSL}(2,7)$, which, as mentioned earlier, satisfies the hypothesis of the theorem. If $q=3$, then according to \cite{Con}, $G$ has a maximal subgroup isomorphic to $S_4$ of index $234=2\cdot 3^2\cdot 13$. If $q=4$, then by \cite{Con}, $G$ has a maximal subgroup isomorphic to $\mathrm{PSL}(2,7)$ of index $120=2^3\cdot 3\cdot 5$. In either case, the index is divisible by three distinct primes, contrary to the hypothesis. Henceforth, we assume that $q\geq 5$.

By \cite[Table 8.3]{BHR}, $G$ has a maximal subgroup $M\cong (C_{q^2+q+1}: C_3)/Z$. Since
\[
|G|=\frac{q^3(q^3-1)(q^2-1)}{|Z|},
\]
we obtain
\[
|G:M|=\frac{q^3(q-1)^2(q+1)}{3}.
\]
Note that $(q,q-1)=(q,q+1)=1$ and $(q+1,q-1)\leq 2$. We distinguish three possibilities.

Suppose first that $3\mid q$. Since $q\geq 5$ is a power of $3$, we have $q=3^f$ with $f\geq 2$, and hence
\[
\pi(|G:M|)=\pi(q(q-1)(q+1)).
\]
Given that $|\pi(|G:M|)|\leq 2$,  and $q-1,q+1>1$, the integers $q-1$ and $q+1$ must be powers of the same prime. Because $(q-1,q+1)\leq 2$, this prime must be $2$, which forces $q=3$, contradicting $q\geq 5$.

Suppose next that $3\mid (q-1)$. Since $q\geq 5$ is coprime to $(q-1)^2(q+1)/3>1$, by hypothesis $(q-1)^2(q+1)/3$ is a prime power. As $3\mid (q-1)$, we have $3\mid (q-1)^2/3$, so both $(q-1)^2(q+1)/3$ and $(q-1)^2(q+1)$ are powers of $3$. Consequently,  $3$ divides both $q-1$ and $q+1$, contradicting $(q+1,q-1)\leq 2$.

Suppose finally that $3\mid (q+1)$. Set $\alpha=(q+1)/3$ and $\beta=q-1$. Since $q\geq 5$, we have $1<\alpha<\beta$. Because $|G:M|=q^3\alpha\beta^2$ has at most two distinct prime divisors and $(q,\alpha\beta)=1$, both $\alpha$ and $\beta$ must be powers of the same prime. Since $(\alpha,\beta)\mid 2$, this common prime is $2$, and as $\alpha<\beta$, we obtain $\alpha\mid \beta$. We conclude that $\alpha=(\alpha,\beta)=2$, which yields $q=5$ and $G=\mathrm{PSL}(3,5)$. However, using \cite[Table 8.3]{BHR}, we find that $G$ has a maximal subgroup
\[
W\cong (Z\times\mathrm{SO}(3,5))/Z\cong\mathrm{SO}(3,5).
\]
Since $|W|=120$, we obtain
\[
|G:W|=\frac{2^5\cdot 3\cdot 5^3\cdot 31}{2^3\cdot 3\cdot 5}=2^2\cdot 5^2\cdot 31,
\]
again a contradiction.

\medskip 
(1.3) $n\geq 4$. Then $G=\mathrm{PSL}(n,q)$, where $n\geq 4$ and $q=s^f$ with $s$ a prime. By the lower part of \cite[Table 4.1]{LL}, the possibilities for $(G,K)$ are: $(\mathrm{PSL}(n,s^f), P_1)$, with $P_1$ a maximal parabolic subgroup, where either $n=t^2$ for a prime $t$ and $q=s$ is prime, or $n$ is prime and $f=n^if_0$ with $f_0=1$ or $f_0$ prime; together with $(\mathrm{PSL}(4,2^{2^i}), P_1)$ for $i\geq 0$. 

\medskip
Assume first that $f>1$. Let $r$ be a prime divisor of $f$ and set $q_0=s^{f/r}$, so that $q=q_0^r$. By \cite[Proposition 4.5.3]{KL} and \cite{BHR}, we consider a maximal subgroup $M= (\mathrm{SL}(n,q_0).\left(\frac{q-1}{q_0-1},n\right))/Z$ of $G$ and compute
\[
|G:M|=\frac{q_0^{n(n-1)(r-1)/2}}{\left(\frac{q-1}{q_0-1},n\right)}\prod_{i=2}^{n}\frac{q^i-1}{q_0^i-1}.
\]
In particular, $q_0 \mid |G : M|$. Setting $h = (n, \frac{q-1}{q_0-1})$, the integer
\[
N = \frac{1}{h} \prod_{i=2}^{n} \frac{q^i - 1}{q_0^i - 1}
\]
is coprime to $q_0$, so the fact  that $|\pi(|G : M|)| \le 2$ implies that $N$ is a prime power. Since $n \ge 4$ and $r \ge 2$, we have $nr \ge 8$, and then Lemma 2.3 ensures that $q_0^{nr} - 1 = q^n - 1$ has a primitive prime divisor $p_1 \ge nr + 1 > n \ge h$. Because $p_1 \nmid q_0^n - 1$ and $p_1 > h$, we obtain $p_1 \mid (q^n - 1)/(q_0^n - 1)$, whence $N = p_1^a$ for some positive integer $a$. On the other hand, since $p_1$ is a primitive prime divisor of $q_0^{nr} - 1$, we have $p_1 \nmid q_0^{ir} - 1 = q^i - 1$ for all $1 \le i \le n - 1$. Thus $\prod_{i=2}^{n-1} \frac{q^i - 1}{q_0^i - 1}$ is coprime to $p_1$, and hence divides $h$. In particular,
\[
\frac{q^{n-1} - 1}{q_0^{n-1} - 1} \le h \le \frac{q - 1}{q_0 - 1},
\]
which is impossible since $n - 1 \ge 3$ and $q > q_0 > 1$. Therefore, $f = 1$ and $q = s$ is a prime. Accordingly, either $n = t^2$ with $t$ prime or $n$ is prime, so $n$ is odd or $n = 4$. We next discuss these cases.

\medskip
(a) Suppose that $n$ is odd. Then $n\geq 5$. Assume first that $q$ is odd. By \cite[Proposition 4.8.4]{KL} and \cite{BHR}, $G$ has a maximal subgroup $Y=(Z\times \mathrm{SO}(n,q))/Z\cong \mathrm{SO}(n,q)$. Then
\[
|G:Y|=\frac{\frac{1}{d}q^{n(n-1)/2}\prod_{i=2}^n(q^i-1)}{q^{((n-1)/2)^2}\prod_{j=1}^{(n-1)/2}(q^{2j}-1)}=\frac{1}{d}q^{\frac{n^2-1}{4}}\prod_{j=1}^{\frac{n-1}{2}}(q^{2j+1}-1).
\]
Since $(q,\prod_{j=1}^{(n-1)/2}(q^{2j+1}-1))=1$, the integer $\frac{1}{d}\prod_{j=1}^{(n-1)/2}(q^{2j+1}-1)$ is a prime power. Since $n\geq 5$, $(q-1)^2$ divides $\prod_{j=1}^{(n-1)/2}(q^{2j+1}-1)$, and as $d\mid q-1$, we have $q-1\mid \frac{(q-1)^2}{d}$. 
Because $q-1> 1$ divides the prime power $\frac{1}{d}\prod_{j=1}^{(n-1)/2}(q^{2j+1}-1)$, both $q-1$ and its divisor $d$ are powers of the same prime.
Thus $\prod_{j=1}^{(n-1)/2}(q^{2j+1}-1)$ is also a prime power. However, $(q^3-1)(q^5-1)$ divides this product, and since $q\geq 3$, Lemma \ref{ppd} ensures that $q^3-1$ and $q^5-1$ admit distinct primitive prime divisors, a contradiction. Consequently, $q=2$.

Assume next that $n$ is a prime. By \cite[Proposition 4.3.6]{KL} and \cite{BHR}, $G$ has a maximal subgroup $T\cong C_{2^n-1}:C_n$ of index
\[
|G:T|=\frac{2^{n(n-1)/2}}{n}\prod_{i=2}^{n-1}(2^i-1).
\]
If $n=5$, then $|G:T|=\frac{2^{10}}{5}(2^4-1)(2^3-1)(2^2-1)=2^{10}\cdot 3^2\cdot 7$, which has three distinct prime divisors, a contradiction. If $n\geq 7$, then $n-1\geq 6$, so $(2^5-1)(2^4-1)(2^3-1)(2^2-1)=3^2\cdot 5\cdot 7\cdot 31$ divides $\prod_{i=2}^{n-1}(2^i-1)$. Dividing by the single prime $n$ still leaves at least three distinct odd prime divisors in $|G:T|$ (alongside $2$), a contradiction.

Assume now that $n=t^2$, where $t$ is an odd prime. By \cite[Proposition 4.7.3]{KL} and \cite{BHR}, $G$ has a maximal subgroup $W\cong \mathrm{PSL}(t,2)^2.C_2$ of index
\[
|G:W|=2^{(t^4-3t^2+2t-2)/2}\frac{\prod_{i=2}^n(2^i-1)}{\prod_{j=2}^t(2^j-1)^2}.
\]
Since  $2$ divides $|G:W|$ and $\frac{\prod_{i=2}^n(2^i-1)}{\prod_{j=2}^t(2^j-1)^2}$ is odd and greater than $1$, it must be a prime power. Notice that $n\geq 9$. Lemma \ref{ppd} then shows that $2^n-1$ and $2^{n-1}-1$ have primitive prime divisors $r_1\geq n+1$ and $r_2\geq n$, respectively, with $r_1\neq r_2$. Since $r_1,r_2\nmid \prod_{j=2}^t(2^j-1)^2$, both $r_1$ and $r_2$ divide $\frac{\prod_{i=2}^n(2^i-1)}{\prod_{j=2}^t(2^j-1)^2}$, which is a contradiction.

\medskip
(b) Suppose that $n=4$. We first dispense with the cases $q=2,3$. We know from  \cite{Con} (and also by Theorem \ref{thm2.1}) that $G={\rm PSL}(4,2)\cong A_8$ satisfies the hypothesis, whereas $G={\rm PSL}(4,3)$ has a maximal subgroup $S_6$ with index $2^3\cdot 3^4\cdot 13$, thereby ruling it out. 

Therefore,  we assume $q\geq4$ in what follows. 
   Table 8.8 of \cite{BHR} asserts that $G$ has a maximal subgroup $W=({\rm SL}(2,q^2):C_{q+1}.C_2)/Z$ 
   of index 
   $$|G:W|=\frac{1}{2}q^4(q^3-1)(q-1).$$
 Now, since $|\pi(|G:W|)|\leq2$ and $q$ is an odd prime, we have that $\frac{1}{2}(q^3-1)(q-1)$ is a prime power.
 Also, both $q^3-1$ and $q-1$ are even, so $\frac{1}{2}(q^3-1)(q-1)$ is even  as well, and hence must be a power of $2$. Likewise,  $(q^3-1)(q-1)=(q-1)^2(q^2+q+1)$ is a power of $2$, which is impossible since $q^2+q+1>1$ is odd.

 \medskip

 (2) $G$ is a symplectic group.  
In view of \cite[Table 4.2]{LL}, we only need to consider the cases: $G={\rm PSp}(4,q)$ with  $q=s^f$, where either $q$ is even, or $q^2=2r^a-1$ for some prime $r$ and some positive integer $a$, or $s=2^r-1$ for some prime $r$; and  $G={\rm Sp}(2m, q)$ with $q=2^f$ and $m\geq3$.  
Recall that
$$
|{\rm PSp}(2m,q)|=\frac{q^{m^2}}{(2,q-1)}
\prod_{i=1}^{m}(q^{2i}-1)
$$
and ${\rm Sp}(2,q)\cong{\rm SL}(2,q)$.

\medskip
(2.1) $ G={\rm PSp}(4,q)$, where either $q$ is even, or $q^2=2r^a-1$ or   $s=2^r-1$   for some prime $r$ and positive integer $a$.  Since $G$ is simple $q\geq 3$ . Let $Z=Z({\rm Sp}(4,q))$

  By \cite[Tables 8.12 and 8.14]{BHR}, $G$ has a maximal subgroup $U\cong ({\rm Sp}(2,q^2):C_2)/Z$. Thus
$$
|G:U|=\frac{q^2}{2}(q-1)(q+1).
$$
We only need to distinguish whether $q$ is odd or even. Suppose first that $q$ is even.  Since $q\geq4$, the integers $\frac{q^2}{2}$, $q-1$, and $q+1$ are all greater than $1$ and are pairwise coprime. This means that $|G:U|$ has at least three distinct prime divisors, a contradiction.

Suppose now that $q$ is odd  and write $|G:U|=q^2\frac{q+1}{2}(q-1)$. Since $(q,q^2-1)=1$, the hypothesis guarantees that $\frac{q+1}{2}  (q-1)$ is a prime power. Thus $\frac{q+1}{2}$ and $q-1$ are powers of the same prime. As $\frac{q+1}{2}\leq q-1$, we have that $\frac{q+1}{2}$ divides $q-1$, which means that $(q+1)\mid 2(q-1)=2(q+1)-4$. We conclude then
that $(q+1)\mid 4$, whence $q=3$ and $G={\rm PSp}(4,3)$.  According to \cite{Con}, the distinct  indices of maximal subgroups of $G$ are $3^3, 2^2\cdot 3^2, 2^3\cdot 5$ and $3^2\cdot5$. Therefore, the group ${\rm  PSp}(4,3)$ satisfies the hypothesis of the theorem. 

 \medskip
 (2.2)  $G={\rm Sp}(2m,q)$, where $q=2^f$ and $m\geq 3$.  Notice that in this case $Z({\rm Sp} (2m,q))= 1$. 
  By \cite[Proposition 4.1.19]{KL} and \cite{BHR}, we know that 
$G$ has a maximal subgroup $$U=[q^a]: ({\rm GL}(t,q)\times {\rm Sp}(2m-2t,q)),$$ where $a=\frac{t}{2}-\frac{3t^2}{2}+2tm$. 
Taking $t=m$, we have
$$
U\cong (C_2)^{fm(m+1)/2}:{\rm GL}(m, q)
\quad {\rm and}\quad
|G:U|=\prod_{i=1}^{m}(q^i+1).
$$
By hypothesis, $|G:U|$ has at most two distinct prime divisors. Since $q=2^f$ and $m\geq 3$, we have
$$
(q+1,q^2+1)=(q^2+1,q^3+1)=1 ,\qquad (q+1,q^3+1)=q+1.
$$
As a  consequence,  $(q+1)(q^3+1)=(q+1)^2(q^2-q+1)$ is a prime power. However, since 
$q+1\leq q^2-q+1$, we conclude that $q+1\mid q^2-q+1=(q+1)(q-2)+3$, which forces $q+1=3$, that is, $q=2$.
Now, if $m\geq4$, then $$(2+1)(2^2+1)(2^4+1)=3\cdot5\cdot 17$$ divides $|G:U|$, contradicting the hypothesis. Therefore, the only possibility is $m=3$ and $G={\rm Sp}(6,2)={\rm PSp}(6,2)$. However, using \cite{Con}, we find that $G$ has a maximal subgroup $H\cong{\rm PSU}(3, 3):C_2$
with $|G:H|=120=2^3\cdot3\cdot5$,
 contradicting the hypothesis as well.

\medskip

(3) $G={\rm PSU}(n,q)$ is a unitary group, where $n\geq3$ and $q=s^f$ with $s$ a prime. Let $Z=Z({\rm SU}(n,q))$. 
Note that $$|G|=\frac{1}{d}q^{n(n-1)/2}\prod_{i=2}^n(q^i-(-1)^i), $$ where $d=(q+1,n)$.

\medskip
(3.1) $n=3$. In this case $q\geq3$ since ${\rm PSU}(3,2)$ is not simple.  
 If $q=3$ or $4$, then by \cite{Con} the distinct maximal subgroup indices of ${\rm PSU}(3,3)$  and ${\rm PSU}(3,4)$ are $28,36,63$ and $65,208,416,1600$, respectively. Thus both groups satisfy the hypothesis. If $q=5$, then  $G={\rm PSU}(3,5)$. According to \cite{Con}, $G$ has a maximal subgroup with index $126=2\cdot3^2\cdot7$, a contradiction.

\medskip
Assume henceforth that $q\geq7$. By \cite[Table 8.5]{BHR}, we know that $G$ has a maximal subgroup $W=(C_{q^2-q+1}:C_3)/Z$, and then
$$
|G:W|=\frac{q^3(q-1)(q+1)^2}{3}.
$$
Since $s\mid |G:W|$, by hypothesis we have that the $s'$-part of $|G:W|$, that is, $(q-1)(q+1)^2 /3$ is a prime power. We distinguish three cases.

If $3\mid q$, then $(q-1)(q+1)^2$ is a prime power. As $1<q-1<q+1$ and both are powers of the same prime, we have $q-1\mid q+1$. Hence $q-1\mid 2$, contrary to $q\geq7$.

If $3\mid q+1$, then $(q-1)(q+1)^2/3$ is a prime power. Since $3\mid(q+1)^2/3$, it must be a power of $3$. This contradicts $q-1>1$ and $3\nmid q-1$.

Finally, suppose that \(3\mid q-1\). Since \((q-1)(q+1)^2/3\) is a prime power and $1<(q-1)/3<(q+1)^2\), we have \((q-1)/3\mid(q+1)^2\). The identity
$ (q+1)^2=(q+3)(q-1)+4
$
gives $(q-1)/3\mid4$, so $q\in\{7,13\}$ (since $q\geq 7$). Now, if $q=13$, then $(q-1)(q+1)^2/3=2^4\cdot7^2$ has two distinct prime divisors, a contradiction. On the other hand, if $q=7$, then by \cite{Con}, we know that $G$ has a maximal subgroup of order $336$ and index $16856=2^3\cdot7^2\cdot43$, again contradicting the hypothesis.

\medskip

(3.2) $n=4$. Then $G={\rm PSU}(4,q)$. If $q = 2$, then $G = {\rm PSU}(4,2)  \cong {\rm PSp}(4,3)$, which satisfies the hypothesis of the theorem as shown in (2). Suppose then that  $q\geq 3$. In this case, $G$ possesses a maximal subgroup $W={\rm GU}(3,q)/Z$ by \cite[Table 8.10]{BHR}. It follows that
$$|G:W|=q^3(q^2+1)(q-1).$$ 
The fact that $q$ is coprime to $(q^2+1)(q-1)$ implies that $(q^2+1)(q-1)$ is a prime power. As $1<q-1<q^2+1$, we have $q-1\mid q^2+1=(q+1)(q-1)+2$. Hence $q=3$ and $G={\rm PSU}(4,3)$. But $G$ has a maximal subgroup isomorphic to ${\rm PSU}(4,2)$ with index $126=2\cdot 3^2\cdot 7$ (\cite{Con}), a contradiction.
\medskip 

To analyze the remaining cases for $n$, we proceed as follows. Again using \cite{Con}, we can rule out the case $G={\rm PSU}(6,2)$ (for instance, it has a maximal subgroup isomorphic to  ${\rm PSU}(5,2)$ of index $672 = 2^5 \cdot 3 \cdot 7)$. Henceforth, by \cite[Table 4.3]{LL}, we only need to  consider the case where $n\geq 5$  and $n$ is a prime.

\medskip
(3.3) $n\geq5$ is a prime. 
If $G = {\rm PSU}(5,2)$, then by \cite{Con}, $G$ has a maximal subgroup $S_3 \times 3_+^{1+2}:2A_4$ of index $3520 = 2^6 \cdot 5 \cdot 11$, which excludes this case. Assume henceforth that $(n,q) \neq (5,2)$.

By \cite[Proposition 4.3.6]{KL} and \cite{BHR},  $G$ has a maximal subgroup $W=(\frac{q^n+1}{q+1}:C_n)/Z$. Then we easily compute
$$|G:W|=\frac{1}{n}q^{n(n-1)/2}\prod_{i=1}^{n-1}(q^i-(-1)^i).$$
If $n\mid q^{n(n-1)/2}$, then the hypothesis requires $\prod_{i=1}^{n-1}(q^i-(-1)^i)$ to be a prime power. In particular, we have that $q^4-1$ is a prime power. Thus $(q^2-1)\mid (q^2+1)=(q^2-1)+2$. which gives $q^2=3$, a contradiction. Therefore, $n$ must  divide  $\prod_{i=1}^{n-1}(q^i-(-1)^i)$, so $\frac{1}{n} \prod_{i=1}^{n-1}(q^i-(-1)^i)$ is also a prime power.  If $n\mid q^2-1$, then $(q^4-1)(q^3+1)((q^2-1)/n)$ is a prime power. In this case, $q^4-1$ is a prime power as well, so $q^2-1\mid q^2+1$, forcing $q^2=3$, which is impossible. Accordingly, $n\nmid q^2-1$, and this guarantees that $q^2-1$ is a prime power. As a consequence, $(q-1)\mid (q+1)=(q-1)+2$, whence $q=2$ or $3$, which we analyze next.

\medskip
In the case $q=3$, since $n\geq 5$, we have that $(3^4-1)(3^3+1)(3^2-1)(3+1)=2^{11}\cdot 5\cdot 7$ divides $\prod_{i=1}^{n-1}(q^i-(-1)^i)$. After dividing by the single prime $n$, the factor $\frac{1}{n}\prod_{i=1}^{n-1}(q^i-(-1)^i)$ is still divisible by at least two distinct primes, so we get a contradiction. It only  remains to consider $q=2$; since $(n,q)\neq (5,2)$ and $n$ is a prime, we have $n\geq 7$. Then $(2^6-1)(2^5+1)(2^4-1)(2^3+1)(2^2-1)(2+1)=3^8\cdot 5\cdot 7\cdot 11$ divides $\prod_{i=1}^{n-1}(q^i-(-1)^i)$, so  $\frac{1}{n}\prod_{i=1}^{n-1}(q^i-(-1)^i)$ has at least three distinct prime divisors, again a contradiction.

\medskip
(4) $G$ is an orthogonal group.

\medskip
(4.1) $G={\rm P\Omega}_{2m+1}(q)$, where $m\geq3$ and $q=s^f$ with $s$ an odd prime.
We know that $$|G|= \frac{q^{m^2}}{d}\prod_{i=1}^m(q^{2i}-1),$$ where $d=(2,q-1)=2$.  In addition, in this case  $Z({\rm \Omega}_{2m+1}(q))=1$.

By \cite[Proposition 4.1.20]{KL} and \cite{BHR}, we have that $G$ possesses a maximal subgroup $W=[q^a]:\frac{1}{2}{\rm GL}_m(q)$, where $a=m(2m+1)-\frac{m}{2}(3m+1)=\frac{1}{2}m(m+1)$. Then 
$$|G:W|= \prod_{i=1}^m(q^i+1).$$
As $m \ge 3$, the product $(q+1)(q^2+1)(q^3+1) = (q+1)^2(q^2+1)(q^2-q+1)$ divides $|G:W|$. It is clear that $(q^2+1, q^2-q+1) = (q, q^2-q+1) = 1$. Furthermore, because $q \ge 3$ is odd, $q^2-q+1 > 1$ is odd (and thus has at least one odd prime divisor), whereas $q^2+1$ is even. Taking into account that $|\pi(|G:W|)| \le 2$, it follows that $q^2+1$ must be a power of $2$. However, the fact that $q^2+1 \equiv 2 \pmod 4$ forces $q^2+1 = 2$, clearly a contradiction.

\medskip
(4.2) $G={\rm P\Omega}^\pm_{2m}(q)$, where $m\geq4$ and $q=s^f$ with $s$ prime. Let $Z=Z({\rm \Omega}^\pm_{2m}(q))$.

Note that $$|G|=\frac{1}{d}q^{m(m-1)}(q^m\mp 1)\prod_{i=1}^{m-1}(q^{2i}-1),$$ where $d=(4,q^m\mp 1)$. Let $e=(2,q-1)$. Observe also that $|Z|=\frac{d}{e}$.

\medskip
(a) Assume first that $q$ is odd. Then $e=2$.
By \cite[Proposition 4.1.6]{KL} and \cite{BHR}, we know that $G$ has a maximal subgroup $W=({\rm \Omega}_{2m-1}(q).C_2)/Z$. Using $$|{\rm \Omega}_{2m-1}(q)|=\frac{1}{2}q^{(m-1)^2}\prod_{i=1}^{m-1}(q^{2i}-1),$$ we  obtain $$|G:W|=\frac{1}{2} q^{m-1}(q^m\mp1).$$ Since $q$ is coprime to $(q^m\mp1)/2$, it follows that $(q^m\mp1)/2$ is a prime power. If $4\mid (q^m\mp1)$, then   $q^m\mp1=2^t$, that is, $s^{fm}\mp1=2^t$, where $t$ is a positive integer. Since $q\geq3$, we have $t\geq2$. If $s^{fm}-1=2^t$, then, by Lemma \ref{ppd}, we have that $s^{fm}-1 $  has a primitive prime divisor $r$ with $r\geq fm+1$, a contradiction. If $s^{fm}+1=2^t$, then, again by Lemma \ref{ppd}, $s^{2fm}-1$ has a primitive prime divisor $r_0$ with $r_0\geq 2fm+1$, which implies that $r_0\mid s^{fm}+1$, also a contradiction. Accordingly, we can assume  that    either $G={\rm P\Omega}^-_{2m}(q)$ and  $4\nmid (q^m +1)$, or $G={\rm P\Omega}^+_{2m}(q)$ and $4\nmid q^m-1$.

\medskip
Consider first the case where $G={\rm P\Omega}^-_{2m}(q)$ and  $4\nmid (q^m +1)$. Then $d=(4,q^m+1)=2=e$ and $Z=1$. By hypothesis, $(q^m+1)/2$ is an odd prime power, so we can write $q^m+1=2r^t$, with $r$  an odd prime and $t$ a positive integer. Suppose that $m$ is not a  power of $2$, and write $m=m_1m_2$, where $m_1$ is an odd prime  and $m_2\geq 1$. Then we have $(q^{m_2}+1)\mid (q^m+1)= 2r^t$ and thus $r\mid  (q^{m_2}+1)$. Now, since $s$ is odd and $2fm\geq 8$, Lemma \ref{ppd} shows that $s^{2fm}-1$ has a primitive prime divisor $p_0$
 with $p_0\geq 2fm+1\geq 9$.  From $s^{2fm}-1=(s^{fm}+1)(s^{fm}-1)$, we obtain that $p_0\mid (s^{fm}+1)=q^m+1$, and hence $p_0=r$. On the other hand, since $r\mid q^{m_2}+1$, we have $r\mid (q^{2m_2}-1)=s^{2fm_2}-1$. Because $2fm_2< 2fm$, this  contradicts the fact that $p_0=r$ is a primitive prime divisor of $s^{2fm}-1$. As a result, $m$ is  a power of $2$.

We now consider another maximal subgroup. By \cite[Proposition 4.3.16]{KL} and \cite{BHR}, we know that $G$ has a maximal subgroup $U = ({\rm P\Omega}_m^-(q^2).C_2)/Z \cong {\rm P\Omega}_m^-(q^2).C_2$. Here,
$$
|U| = \frac{2}{d_1} q^{m(\frac{m}{2}-1)}(q^m+1)\prod_{i=1}^{\frac{m}{2}-1} (q^{4i}-1),
$$
where $d_1 = (4, (q^2)^{\frac{m}{2}}+1)=(4, q^m+1)= d$. Then
$$|G : U| = \frac{1}{2} q^{\frac{m^2}{2}}N, \quad {\rm with }\quad N= \frac{\prod_{i=1}^{m-1} (q^{2i}-1)}{\prod_{j=1}^{\frac{m}{2}-1} (q^{4j}-1)}.
$$
It is clear that $N$ must be a prime power. Since $\frac{1}{2}(q^{2(m-1)}-1)(q^{2(m-3)}-1)$ divides $|G : U|$, we get that $\frac{1}{2}(q^{2(m-1)}-1)(q^{2(m-3)}-1)$ is a prime power. Because $q$ is odd, $q^2 \equiv 1 \pmod{8}$, so both factors $q^{2(m-1)}-1$ and $q^{2(m-3)}-1$ are divisible by $8$; so $\frac{1}{2}(q^{2(m-1)}-1)(q^{2(m-3)}-1)$ is even, and hence a power of $2$,  which means that $q^{2(m-1)}-1$ and $q^{2(m-3)}-1$ are powers of $2$. This implies that $(q^{2(m-3)}-1)$ divides $(q^{2(m-1)}-1)$. Since $(q^{2(m-1)}-1, q^{2(m-3)}-1) = q^2 - 1$, we obtain $q^2 - 1 = q^{2(m-3)} - 1$, and so $m = 4$. In particular, it follows that $q^6 - 1$ is a power of $2$. Now, as $q^6 - 1 = (q^3 - 1)(q^3 + 1)$, we conclude that $(q^3 - 1) \mid (q^3 + 1)$, and hence $q^3 = 3$, a contradiction.

 \medskip
Assume now that $G={\rm P\Omega}^+_{2m}(q)$ and  $4\nmid (q^m -1)$. As $q^2\equiv 1$  $({\rm mod}$ $8)$,  $m$ is necessarily odd (so $m\geq 5$). Also, $d=(4,q^m-1)=2=e$, so $Z=1$.  By our assumption, $(q^m-1)/2$ is an odd prime power, so  we write  $q^m-1=2p^a$, where $p$ is an odd prime and $a$ is a positive integer. Since $s$ is odd and  $fm\geq 5$, there exists a primitive prime divisor  $p_0\geq fm+1$ of  $s^{fm}-1$ by  Lemma \ref{ppd}, so necessarily  $p_0=p$.  Since $(s-1)\mid (s^{fm}-1)=2p^a$ and $p_0\nmid (s-1)$, we deduce that $s-1=2$, whence $s=3$. If $f>1$, then $(3^f-1)\mid (3^{fm}-1)=2p^a$; but since $fm> f$, we have $3^f-1=2$, a contradiction. Therefore, $f=1$ and $q=3$. By the same argument, if $m$ were composite, any prime divisor $m_1$ of $m$ would satisfy
$1< m_1< m$ and give $3^{m_1}- 1= 2$, a contradiction. Thus $m\geq  5$ is an odd prime.

Now, by using \cite[Proposition 4.3.20]{KL} and \cite{BHR}, we find that $G$ has a maximal subgroup $T=({\rm\Omega}_{m}(q^2).C_2)/Z\cong {\rm\Omega}_{m}(q^2).C_2$, whose order is $$|T|=2|{\rm\Omega}_{m}(q^2)|= \frac{2q^{\frac{(m-1)^2}{2}}}{d_2}\prod_{i=1}^{\frac{m-1}{2}}(q^{4i}-1),$$ where $d_2=(2,q^2-1)=2$. Hence, $$|G:T|=\frac{1}{2}q^{\frac{m^2-1}{2}}(q^m-1)\frac{\prod_{i=1}^{m-1}(q^{2i}-1)}{\prod_{j=1}^{\frac{m-1}{2}}(q^{4j}-1)}=p^a q^{\frac{m^2-1}{2}} \frac{\prod_{i=1}^{m-1}(q^{2i}-1)}{\prod_{j=1}^{\frac{m-1}{2}}(q^{4j}-1)}.$$
 Since $q=3$ is coprime to the remaining factors of $|G:T|$, by  our assumptions, the factor $$\frac{\prod_{i=1}^{m-1}(3^{2i}-1)}{\prod_{j=1}^{\frac{m-1}{2}}(3^{4j}-1)}$$ must be a power of $p$. This leads to a contradiction, since for $m\geq 5$  this factor  is divisible by $(3^2-1)(3^6-1)=2^6\cdot 7\cdot 13$.

\medskip
(b) Suppose now that $q$ is even.  
 By \cite[Proposition 4.1.20]{KL} and [1], $G$ has a maximal subgroup 
 $J = [q^{a}] : (GL_{2}(q) \times \Omega_{2m-4}^{\pm}(q))$, where $a = 4m - 7$. Since $q$ is even, we have $d = e = 1$ and $Z = 1$. Recall that
$$
|\Omega_{2m-4}^{\pm}(q)| = q^{(m-2)(m-3)}(q^{m-2} \mp 1)\prod_{i=1}^{m-3}(q^{2i} - 1),
$$
and hence
$$
|G : J| = \frac{(q^{m} \mp 1)(q^{2(m-1)} - 1)(q^{m-2} \pm 1)}{(q - 1)(q^{2} - 1)}.
$$
By hypothesis, $|G : J|$ has at most two distinct prime divisors. When $q = 2$ and $m \in \{4, 5, 6, 7\}$, a direct computation shows that $|G : J|$ is divisible by at least three distinct primes. Thus, we may assume henceforth that $fm \ge 8$.

\medskip
Assume first that $G = P\Omega_{2m}^{+}(q)$, so
$$
|G : J| = \frac{(q^{m} - 1)(q^{2(m-1)} - 1)(q^{m-2} + 1)}{(q - 1)(q^{2} - 1)}.
$$
If $m = 4$, then $|G : J| = (q^{2} + 1)^{2}(q^{3} + 1)(q^{2} + q + 1)$. Since $q$ is even, the three integers $q^{2} + 1$, $q^{3} + 1$, and $q^{2} + q + 1$ are greater than $1$ and pairwise coprime, so $|G : J|$ has at least three distinct prime divisors, a contradiction. Hence $m \ge 5$. Since $fm \ge 8$ and $m \ge 5$, each of the exponents $fm$, $2f(m-1)$, and $2f(m-2)$ is at least $8$ and strictly greater than $2f$. By Lemma 2.3, $2^{fm} - 1 = q^{m} - 1$, $2^{2f(m-1)} - 1 = q^{2(m-1)} - 1$, and $2^{2f(m-2)} - 1 = q^{2(m-2)} - 1$ possess primitive prime divisors $r_{1}$, $r_{2}$, and $r_{3}$, respectively. Because the orders $fm$, $2f(m-1)$, and $2f(m-2)$ of $2$ modulo $r_{1}$, $r_{2}$, and $r_{3}$ are distinct and greater than $2f$, these primes are distinct and coprime to $(q-1)(q^{2}-1)$. Since $r_{3} \nmid (q^{m-2}-1)$, we also have $r_{3} \mid (q^{m-2}+1)$, so $r_{1}, r_{2}, r_{3}$ divide $|G : J|$, a contradiction

\medskip
Assume now that $G = P\Omega_{2m}^{-}(q)$, so
$$
|G : J| = \frac{(q^{m} + 1)(q^{2(m-1)} - 1)(q^{m-2} - 1)}{(q - 1)(q^{2} - 1)}.
$$
If $q = 4$ and $m = 4$, then $|G : J| = (4^{4} + 1)(4^{6} - 1)/(4 - 1) = 3 \cdot 5 \cdot 7 \cdot 13 \cdot 257$, contradicting the hypothesis. Since $m \ge 4$ and $fm \ge 8$, we therefore have $(f, m) \ne (2, 4)$, so that $f(m - 1) \ge 3f \ge 3$ and $f(m - 1) \ne 6$. By Lemma 2.3, $2^{2fm} - 1$, $2^{2f(m-1)} - 1$, and $2^{f(m-1)} - 1$ admit primitive prime divisors $r_{1}$, $r_{2}$, and $r_{3}$, respectively. 
Since the orders $2fm$, $2f(m-1)$, and $f(m-1)$ of $2$ modulo $r_{1}$, $r_{2}$, and $r_{3}$ are distinct and greater than $2f$, these primes are distinct and coprime to $(q-1)(q^{2}-1)$. Moreover, $r_{1} \nmid (q^{m}-1)$ yields $r_{1} \mid (q^{m}+1)$, while $r_{2}, r_{3} \mid (q^{2(m-1)}-1)$. Thus, $r_{1}$, $r_{2}$, and $r_{3}$ are distinct prime divisors of $|G : J|$, a contradiction.

\medskip
(5) $G$ is an exceptional group of Lie type.  By \cite[Table 5.1]{LL}, we only need to consider the groups $^2{\rm F}_4(2)', {\rm G}_2(3), {\rm G}_2(4)$, or $^2{\rm B}_2(q)$, where $q=2^f$ with $f$ an odd prime. By \cite{Con}, we can rule out the first three groups. For instance, ${}^2{\rm F}_4(2)'$ has a maximal subgroup $2^2.[2^8]:S_3$ of index $2925=3^2\cdot 5^2\cdot 13$; ${\rm G}_2(3)$ has a maximal subgroup $(3_+^{1+2}\times 3^2):2S_4$ of index $364=2^2\cdot 7\cdot 13$; and ${\rm G}_2(4)$ has a maximal subgroup ${\rm PSU}(3,4):2$ of index $2016=2^5\cdot 3^2\cdot 7$. Suppose then that $G=^2{\rm B}_2(q)$ and recall that $$|^2{\rm B}_2(q)|=q^2(q^2+1)(q-1).$$
By \cite[Theorem 4.1]{Wil}, $G$ has a maximal subgroup $W=C_{q+\sqrt{2q}+1}:C_4$. Then $$|G:W|=\frac{q^2}{4}(q-1)(q-\sqrt{2q}+1).$$

Since $q=2^f$ with $f\geq 3$, it is clear that the three integers $\frac{q^2}{4}$, $q-1$, and $q-\sqrt{2q}+1$ are all strictly greater than $1$. Moreover, $\frac{q^2}{4}$ is a power of $2$ whereas both $q-1$ and $q-\sqrt{2q}+1$ are odd (as $q-\sqrt{2q}+1$ divides $q^2+1$), and $(q^2+1,q-1)=((q-1)(q+1)+2,q-1)=(2,q-1)=1$, which implies $(q-\sqrt{2q}+1,q-1)=1$. We conclude that $\frac{q^2}{4}$, $q-1$, and $q-\sqrt{2q}+1$ are pairwise coprime integers greater than $1$, so $|G:W|$ is divisible by at least three distinct primes, a contradiction.
\end{proof}

We can now give the proof of Theorem A.

\begin{proof}[Proof of Theorem A]
The theorem follows from Theorems \ref{thm2.1}, \ref{thm2.2} and \ref{thm2.4}, taking into account the isomorphisms ${\rm PSL}(2,5)\cong A_5$, ${\rm PSL}(2,9)\cong A_6$ and ${\rm PSL}(4,2)\cong A_8$.
\end{proof}
\bigskip
\noindent
{\bf Acknowledgments}.
This work was supported by the National Natural Science Fund of China (No.  12071181). The second author gratefully acknowledges the warm hospitality received during his research stay at the Universitat Jaume I, Spain, in autumn 2026. 

\bigskip
\noindent 
{\bf Statements and declarations}.
The authors have no competing interests to declare that are relevant to
the content of this article. Data sharing is not applicable to this article, as no datasets were used or generated during this
study.

\bibliographystyle{plain}

\end{document}